\documentclass[12pt]{amsart}
\input{diagrams}

\newtheorem{proposition}{Proposition}[section]
\newtheorem{lemma}[proposition]{Lemma}
\newtheorem{corollary}[proposition]{Corollary}
\newtheorem{theorem}[proposition]{Theorem}

\theoremstyle{definition}

\theoremstyle{remark}

\newtheorem{remarks}[proposition]{Remarks}

\newcommand{\thlabel}[1]{\label{th:#1}}
\newcommand{\thref}[1]{Theorem~\ref{th:#1}}
\newcommand{\selabel}[1]{\label{se:#1}}
\newcommand{\seref}[1]{Section~\ref{se:#1}}
\newcommand{\lelabel}[1]{\label{le:#1}}
\newcommand{\leref}[1]{Lemma~\ref{le:#1}}
\newcommand{\prlabel}[1]{\label{pr:#1}}
\newcommand{\prref}[1]{Proposition~\ref{pr:#1}}
\newcommand{\colabel}[1]{\label{co:#1}}
\newcommand{\coref}[1]{Corollary~\ref{co:#1}}
\newcommand{\relabel}[1]{\label{re:#1}}

\newcommand{\Hom}{\rm{Hom}}

\newcommand{\Ext}{\rm{Ext}}

\def\ot{\otimes}

\newcommand{\Mm}{\mathcal{M}}

\def\text#1{{\rm {\rm #1}}}

\begin{document}
\title[Projectivity of a relative Hopf module]{Projectivity of a relative Hopf module
over the subring of coinvariants}
\author{S. Caenepeel}
\address{Faculty of Applied Sciences,
Vrije Universiteit Brussel, VUB, B-1050 Brussels, Belgium}
\email{scaenepe@vub.ac.be}
\urladdr{http://homepages.vub.ac.be/~scaenepe/}
\author{T. Gu\'ed\'enon}
\address{Faculty of Applied Sciences,
Vrije Universiteit Brussel, VUB, B-1050 Brussels, Belgium}
\email{guedenon@caramail.com}

\begin{abstract}
Let $k$ be a commutative ring, $H$ a faithfully flat Hopf algebra with bijective
antipode, $A$ a $k$-flat right $H$-comodule algebra. We investigate when a relative
Hopf module is projective over the subring of coinvariants $B=A^{{\rm co}H}$, and we
study the semisimplicity of the category of relative Hopf modules. 
\end{abstract}
\maketitle

\section*{Introduction}\selabel{0}
Let $H$ be a Hopf algebra with bijective antipode, and $A$ a left
$H$-module algebra. We can then form the smash product $A\# H$.
In \cite{Thomas}, the second author investigated necessary and
sufficient conditions for the projectivity of an $A\# H$-module $M$
over the subring of invariants $A^H$.\\
Our starting point is the following: if $H$ is finitely
generated and projective, then $H^*$ is also a Hopf algebra, $A$ is
a right $H^*$-comodule algebra, and $A^H=A^{{\rm co}H^*}$. The
category of left $A\# H$-modules is isomorphic to ${}_A\Mm^{H^*}$, the
category of relative $(A,H^*)$-Hopf modules, i.e. $k$-modules together with
a left $A$-action and a right $H^*$-coaction, satisfying an appropriate
compatibility condition. Thus \cite{Thomas} brings us necessary
and sufficient conditions for a relative Hopf module to be projective
as a module over the ring of coinvariants. In this paper, we will
generalize these results to
relative $(A,H)$-Hopf modules, where $H$ is an arbitrary Hopf algebra with
bijective antipode over a commutative ring $k$, and $A$ is a 
($k$-flat) right $H$-comodule algebra.\\ 
Our main result is \thref{2.1}, where we give necessary and sufficient
conditions for projectivity of a relative Hopf module over the subring
of coinvariants $B=A^{{\rm co}H}$. The main tool - based on the methods
developed in \cite{Garcia} and \cite{Thomas} - 
is the basic fact that
the canonical structure of $A$ as left-right $(A, H)$-Hopf module is such that
${}_A\Hom^H(A, A)$, the $k$-module consisting of $A$-linear and $H$-colinear maps, is
isomorphic to $B$ (see \seref{1}). The result can be improved if we
assume that there exists a total integral $\phi :\ H\rightarrow A$,
see \prref{2.2}. In \prref{2.5}, we look at relative modules that
are coinvariantly generated, and we
will present some conditions that
are sufficient, but in general not necessary for projectivity. These conditions
have the advantage that they are easier to verify than the ones from \thref{2.1};
they turn out to be necessary if the coinvariants functor is exact.\\
In Sections \ref{se:3} and \ref{se:4}, we will work over a field $k$. Our
methods will be applied to discuss properties of injective and projective
dimension in the category of relative Hopf modules (\seref{3}), and
semisimplicity of the the category of relative Hopf modules.
Our main result is \coref{4.5}, where we give a sufficient condition
for the category of relative Hopf modules to be semisimple.\\
Let us finally mention that our results may be applied in the following
particular cases:
\begin{itemize}
\item[-] $A=H$ with comodule structure map $\Delta$. In this case, $B$ is
isomorphic to $k$ as
$k$-algebra and trivial $H$-comodule.
\item[-] A crossed product $A=R\#_{\sigma}H$ where $R$ is an algebra with
$H$-action and $\sigma$ is an
invertible map in $Hom_k(H\otimes_kH, R)$ [10, Definition 7.1.1] is a right
$H$-comodule
algebra such that $A^{coH}=R$. The comodule map is given by $a\#_{\sigma}h
: \rightarrow
a\#_{\sigma}h_1 \otimes h_2$.
\item[-] $A={\oplus}_{g \in G}A_g$ is a $G$-graded $k$-algebra for a group $G$ and
$H=kG$ is the group
algebra of $G$. Then $A$ is a right $H$-comodule algebra, $A^{{\rm co}H}=A_1$ and
${}_{A}\Mm^H$ is the category of $G$-graded left $A$-modules.
\item[-] $H=k[G]$, the affine coordinate ring of an affine $k$-group scheme $G$
and $A$ is a right
$k[G]$-comodule algebra. A right $k[G]$-comodule is also called a left
$G$-module (see \cite{donkin}, \cite{Jantzen}). If $G$ is
linearly reductive, meaning that every $G$-module is completely reducible, then the
functor $(-)^{{\rm co}k[G]}:\ \Mm^{k[G]}\to
{}_{k}\Mm$ is exact. In the particular case of an affine algebraic group $G$
over an algebraically
closed field $k$ acting rationally on $A$, ${}_{A}\Mm^{k[G]}$ is the
category of rational $(A,G)$-modules (see \cite{Doraiswamy}, \cite{Magid}).
\end{itemize}
For more detail on Hopf algebras and the category of relative Hopf modules,
we refer to the literature, see for example \cite{Abe}, \cite{Caenepeel98},
\cite{CaenepeelMZ02}, \cite{DascalescuNR}, \cite{Montgomery}, \cite{Sweedler}. 

\section{Preliminary Results}\selabel{1}
Troughout this paper, $k$ will be a commutative ring, and $H$ will
be a Hopf algebra with bijective antipode $S$. We freely use the Sweedler-Heyneman
notation for the comultiplication:
$$\Delta(h)=h_1\ot h_2$$
We use a similar notation for the right $H$-coaction $\rho$ on a right
$H$-comodule $M$:
$$\rho(m)=m_0\ot m_1$$
$\Mm^H$ will be the category of right $H$-comodules and right $H$-colinear
maps. The $k$-module consisting of all right $H$-comodule maps between
two right $H$-comodules $M$ and $N$ will be denoted by $\Hom^H(M,N)$.\\
The tensor product of two right $H$-comodules $M$ and $N$ is again a right
$H$-comodule. The right $H$-coaction on $M\ot N$ is given by
$$\rho(m\ot n)=m_0\ot n_0\ot m_1n_1$$
For any right $H$-comodule $M$,
$$M^{{\rm co}H}=\{m\in M~|~\rho(m)=m\ot 1\}$$
is called the submodule of coinvariants of $M$.\\
A right $H$-comodule algebra $A$ is an algebra in the category $\Mm^H$.
This means that $A$ is a right $H$-comodule and a $k$-algebra such that
the unit and multiplication maps are right $H$-colinear:
$${\rho}_A (ab) = a_0b_0 \otimes a_1b_1 \quad \hbox {and} \quad 
{\rho}_A (1_A) = 1_A \otimes 1_H$$
A (left-right) relative $(A,H)$-Hopf module is a $k$-module $M$ together with
a left $A$-action and a right $H$-coaction such that
$$\rho(am)=a_0m_0\ot a_1m_1$$
for all $a\in A$ and $m\in M$. ${}_A\Mm^H$ is the category of
relative $(A,H)$-Hopf modules and left $A$-linear and right $H$-colinear
maps. For two relative Hopf modules $M$ and $N$, we write
${}_A\Hom^H(M,N)$ for the $k$-module consisting of all 
$A$-linear $H$-colinear maps from $M$ to $N$.\\
$B=A^{{\rm co}H}$ is a $k$-algebra. The coinvariants $M^{{\rm co}H}$
of a relative Hopf module $M$ form a $B$-module.

\begin{lemma}\lelabel{1.1}
\begin{enumerate}
\item Take $M\in \Mm^H$ and $N\in {}_{A}\Mm^H$.
\begin{enumerate}
\item[a)] $N\ot M\in {}_{A}\Mm^H$, with left $A$-action given by
$$a(n\otimes m)=(an)\otimes m$$
\item[b)] If $H$ is commutative, then $M\otimes N\in {}_{A}\Mm^H$, with left $A$-action
given by
$$a(m\otimes n)=m\otimes an$$
\end{enumerate}
\item Let $A$ be commutative, and take $M,N\in {}_{A}\Mm^H$. Then
$M\ot_A N\in {}_{A}\Mm^H$, under the coaction
$${\rho}(m\otimes n)=m_0 \otimes n_0\otimes m_1n_1$$ 
\item Take $M,N\in {}_{A}\Mm^H$, and assume that
$M$ is finitely generated 
projective as a left $A$-module.
\begin{enumerate}
\item[a)] ${}_A\Hom(M, N)\in\Mm^H$ and 
$${}_A\Hom^H(M, N)={}_AHom(M,N)^{{\rm co}H}$$  
\item[b)] $N\cong {}_A\Hom(A, N)\in{}_A\Mm^H$, with left $A$-action
$$(af)(u)=f(ua)$$
\item[c)] If $A$ is commutative, then
${}_A\Hom(M, N)\in{}_A\Mm^H$.
\end{enumerate}
\end{enumerate}
\end{lemma}

\begin{proof}
1a) The only nontrivial thing that we have to show is the compatibility of
the right $H$-coaction with the left $A$-action. For any $a \in A$, $m \in M$ and
$n\in N$, we have
$${\rho}_{N \otimes M}((an)\otimes m)=a_0n_0\otimes m_0\otimes a_1n_1m_1=a_0(n\otimes m)_0\otimes
a_1(n\otimes m)_1$$
1b) As in 1a), we have to prove compatibility. For 
$a \in A$, $m \in M$ and
$n\in N$, we have
\begin{eqnarray*}
&&\hspace*{-2cm}
{\rho}_{M \otimes N}(m\otimes an)=m_0\otimes
a_0n_0\otimes m_1a_1n_1\\
&=& m_0\otimes a_0n_0\otimes a_1m_1n_1 =a_0(m\otimes n)_0\otimes
a_1(m\otimes n)_1
\end{eqnarray*}
2) The right $H$-coaction on $M\ot_AN$ is well-defined since
$$\rho(ma\ot n)=m_0a_0\ot n_0\ot m_1a_1n_1=\rho(m\ot an)$$
and it is straightforward to verify that $M\ot_AN$ is a right
$H$-comodule. Finally
$${\rho}_{M\otimes_AN}(am\otimes n)=a_0m_0 \otimes n_0\otimes a_1m_1n_1=
a_0(m\otimes n)_0 \otimes a_1(m\otimes n)_1$$
3a) The fact that $M$ is finitely generated projective as a left $A$-module
implies that we have a natural isomorphism
$${}_A\Hom(M, N) \otimes H \cong {}_A\Hom(M, N \otimes H)$$
Using this isomorphism, we define a map
$${\pi}:\ {}_A\Hom(M, N)
\rightarrow {}_A\Hom(M, N) \otimes H; \quad \pi(f)= f_0 \otimes f_1$$
by
$$({\pi}(f))(m)=f_0(m) \otimes f_1=f(m_0)_0 \otimes S^{-1}(m_1)f(m_0)_1$$
defining a right $H$-coaction on ${}_A\Hom(M, N)$.\\
Take $f \in {}_A\Hom^H(M, N)$ and $m \in M$. Then
$$f(m)_0\otimes f(m)_1=f(m_{0}) \otimes m_{1}$$ 
so 
\begin{eqnarray*}
&&\hspace*{-2cm}{\pi}(f)(m)=f_0(m) \otimes f_1 = f(m_0)_0 \otimes
S^{-1}(m_1)f(m_0)_1\\
&=&
f(m_0) \otimes S^{-1}(m_2)m_1
=f(m_0) \otimes \varepsilon(m_1)1= (f \otimes 1)(m)
\end{eqnarray*}
proving that $f$ is coinvariant. Conversely, take a coinvariant
$f \in {}_A\Hom(M, N)^{{\rm co}H}$. For every $m \in M$ we have
$$f_0(m)
\otimes f_1 =f(m_0)_0
\otimes S^{-1}(m_1)f(m_0)_1=f(m) \otimes 1$$ 
and we deduce that
\begin{eqnarray*}
&&\hspace*{-2cm} f(m_0)\ot m_1=f(m_0)\ot m_1 1\\
&=&f(m_0)_0
\otimes m_2S^{-1}(m_1)f(m_0)_1\\
&=& f(m)_0\ot f(m)_1
\end{eqnarray*}
so $f$ is $H$-colinear, as needed.\\
3b) By 3a), ${}_A\Hom(A, N)$ is a right $H$-comodule. The map 
$$\psi:\
{}_A\Hom(A, N)\rightarrow N,~~\psi(f)=f(1)$$
is clearly an isomorphism of $A$-modules. $\psi$ is also
right $H$-colinear, since
\begin{eqnarray*}
&&\hspace*{-2cm}
{\psi(f)}_o \otimes {\psi(f)}_1=f(1)_0 \otimes f(1)_1=f(1_0)_0 \otimes
S^{-1}(1)f(1_0)_1\\
&=&f_0(1)\otimes f_1=\psi(f_0) \otimes f_1
\end{eqnarray*}
3c) ${}_A\Hom(M, N)$ is an $A$-module, and, by 3a), it is a right $H$-comodule. 
The $A$-action and $H$-coaction are compatible since
\begin{eqnarray*}
&&\hspace*{-2cm}
((af)_0 \otimes (af)_1)(m)=((af)(m_0))_0
\otimes S^{-1}(m_1)((af)(m_0))_1\\
&=&a_0(f(m_0)_0)
\otimes S^{-1}(m_1)a_1(f(m_0)_1)\\
&=&a_0(f(m_0)_0)
\otimes a_1S^{-1}(m_1)(f(m_0)_1)=a_0(f_0(m)) \otimes a_1f_1\\
&=& (a_0f_0)(m) \otimes
a_1f_1=(a_0f_0 \otimes a_1f_1)(m)
\end{eqnarray*}
and this proves that
${}_A\Hom(M, N)\in {}_{A}\Mm^H$.
\end{proof}

Consider the functor
$$(-)^{{\rm co}H} :\ {}_{A}\Mm ^H \longrightarrow {}_{B}\Mm $$
$(-)^{coH}$ commutes with direct sums, and has a left adjoint
$$T=A\ot_B -:\ {}_{B}\Mm\longrightarrow{}_{A}\Mm ^H$$
The unit and counit of the adjunction are the following: for
$N\in {}_{B}\Mm $ and $M\in {}_{A}\Mm ^H$:
$$u_N:\ N\longrightarrow ( A\otimes {}_{B}N )^{{\rm co}H},~~
u_{N}(n)= 1\otimes n$$
$$c_M:\ A\otimes_{B}M^{{\rm co}H}\longrightarrow M,~~
c_{M}(a\otimes m)=am$$

\begin{lemma}\lelabel{1.2}
\begin{enumerate}
\item The functor $(-)^{{\rm co}H}$ is naturally isomorphic to
$${}_A\Hom^H(A, -) :\ {}_{A}\Mm ^H \longrightarrow {}_{B}\Mm$$
\item $c_{A^{(I)}}$ is an isomorphism for any set $I$.
\item $u_{B^{(I)}}$ is an isomorphism for any set $I$.
\end{enumerate}
\end{lemma}

\begin{proof}
1) By \leref{1.1} (3), ${}_A\Hom^H(A, M)={}_A\Hom(A, M)^{{\rm co}H}$, and 
${}_A\Hom(A,
M)$ is isomorphic to $M$ in ${}_{A}\Mm ^H$. So ${}_A\Hom^H(A, M)$ is $B$-isomorphic to
$M^{{\rm co}H}$.\\
2) It suffices to observe that
$$(A^{(I)})^{{\rm co}H} = (A^{{\rm co}H})^{(I)}=B^{(I)}$$ 
and
$$c_{A^{(I)}} :\ A\otimes_{B}(A^{{\rm co}H})^{(I)} \longrightarrow A^{(I)}$$ 
is the canonical isomorphism.\\
3) From the fact that $u$ and $c$ are the unit and the counit of the
adjunction, we derive that
$$(c_{A^{(I)}})^{{\rm co}H} \circ u_{(A^{(I)})^{{\rm co}H}}=1_{(A^{(I)})^{{\rm co}H}}$$
and it follows that 
$u_{(A^{(I)})^{{\rm co}H}}=u_{(A^{{\rm co}H})^{(I)}}$ is an isomorphism.
\end{proof}

\begin{lemma}\lelabel{1.3}
Take $M\in {}_{A}\Mm^H$. Then 
\begin{enumerate}
\item $M\otimes H\in {}_{A}\Mm^H$ and 
 $(M\otimes H)^{{\rm co}H}\cong M$ as a $B$-module.
\item $u_M$ is an injection of $B$-modules.
\end{enumerate} 
\end{lemma}

\begin{proof}
1) The structure maps on $M\ot H$ are the following:
$$a(m\ot h)=am\ot h~~;~~\rho(m\ot h)=m_0\ot h_1\ot m_1h_2$$
It is clear that the map 
$$f:\
M\rightarrow M\otimes H;~~
f(m)=m_0\otimes S(m_1)$$
is left $B$-linear. Also
$$
\rho(f(m))= m_0\ot S(m_3)\ot m_1S(m_2)=m_0\ot S(m_1)\ot 1$$
so $f(m)\in (M\ot H)^{{\rm co}H}$, and we have a $B$-linear map
$$f:\ M\to (M\ot H)^{{\rm co}H}$$
Its inverse is 
$$g:\ (M\ot H)^{{\rm co}H}\to M;~~g(m\ot h)=\varepsilon(h)m$$
It is easy to see that $g$ is a left inverse of $f$.
Take $\sum_i m_i\ot h_i\in  (M\ot H)^{{\rm co}H}$. Then
$$f(g(\sum_i m_i\ot h_i))=\sum_i\varepsilon(h_i)f(m_i)=
\sum_i\varepsilon(h_i)m_{i0}\ot S(m_{i1})$$
We also have that
$$\sum_i m_{i0}\ot h_{i1}\ot m_{i1} h_{i2}=\sum_i m_i\ot h_i\ot 1$$
applying $S$ to the third tensor factor, and then multiplying the
second and third factor, we obtain
$$\sum_i \varepsilon(h_i)m_{i0}\ot S(m_{i1})=\sum_i m_i\ot h_i$$
and it follows that $g$ is also a right inverse of $f$.\\

2) Set $W=M\otimes H$. By \leref{1.2}, 
$$(c_{W})^{coH}\circ
u_{W^{coH}}=1_{W^{coH}}$$ 
By 1),
$W^{{\rm co}H}\cong M$ as $B$-modules, so $u_M$ is an injection of
$B$-modules.
\end{proof}

\section{Projectivity of relative Hopf modules}\selabel{2}
We keep the notation of \seref{1}. We will now discuss when a relative
Hopf module is a projective object as a $B$-module.

\begin{theorem}\thlabel{2.1}
For a relative Hopf module
$P\in {}_{A}\Mm^H$, the following
conditions are equivalent:
\begin{enumerate}
\item $P$ is projective in ${}_{B}\Mm$;
\item $A\otimes_{B}P$ is isomorphic in  ${}_{A}\Mm^H$ to a
direct summand of a direct sum of copies of $A$, and $u_P$ is surjective
(bijective);
\item There is a direct summand $M$ in  ${}_{A}\Mm^H$ of a direct sum
of copies of $A$ such that $M^{coH}\cong P$ as a $B$-module.
\end{enumerate}
\end{theorem}

\begin{proof} 
1) $\Rightarrow$ 2). Let $p:\ B^{(I)}
\longrightarrow P$ be a split epimorphism in ${}_{B}\Mm$. Then $$1\otimes
p :\ A\otimes_{B}B^{(I)}\longrightarrow A\otimes_{B}P$$ is a
split epimorphism in ${}_{A}{\Mm^H}$ and $A\otimes_{B}B^{(I)}
\cong A^{(I)}$ in ${}_{A}{\Mm^H}$. Furthermore
$$u_{P}\circ
p= (1\otimes p)^{{\rm co}H}\circ u_{B^{(I)}}$$
 and, by \leref{1.2}, $u_{B^{(I)}}$
is an isomorphism. Thus $u_P$ is an epimorphism
and it follows from \leref{1.3} that $u_P$ is 
bijective.\\
2) $\Rightarrow$ 3). If $u_{P}:\ P \longrightarrow (A\otimes_{B}P)^{{\rm co}H}$
is an epimorphism (and therefore an isomorphism by \leref{1.3}) and if $M$ is a
direct summand in  ${}_{A}\Mm^H$ of a direct sum of copies of $A$
isomorphic to
$A\otimes_{B}P$, then $M$ satisfies the required condition.\\
3) $\Rightarrow$ 1) Let $M$ be a direct summand in  ${}_{A}\Mm^H$
of a direct sum of copies of $A$, such that $M^{{\rm co}H}\cong P$ as
a $B$-module. Then there exists a split epimorphism $f :\
A^{(I)}\longrightarrow M$ in ${}_{A}\Mm^H$. Thus 
$$f^{{\rm co}H} :\
(A^{(I)})^{{\rm co}H}\longrightarrow M^{{\rm co}H}\cong P$$ is a split epimorphism in
${}_{B}\Mm$ and $(A^{(I)})^{{\rm co}H}\cong B^{(I)}$. So $P$ is projective in
${}_{B}\Mm$.
\end{proof}

Recall that a total integral is an $H$-colinear map $\phi:\ H\to A$
such that $\phi(1_H)=1_A$. If there exists a total integral, then 
$u_N$ is an isomorphism of $B$-modules, for every $N\in {}_B\Mm$
(see e.g. \cite[Lemma 23]{CaenepeelMZ02}). From \thref{2.1}, we easily
obtain the following result:

\begin{proposition}\prlabel{2.2}
Assume that there exists a total integral $H\rightarrow A$, and take
$P\in {}_{A}\Mm^H$. $P$ is projective as a $B$-module if and only if
$A\otimes_{B}P$ is
isomorphic in ${}_{A}\Mm^H$ to a  direct summand of a direct sum of
copies of $A$.
\end{proposition}

An $H$-ideal $I$ of $A$ is an $H$-subcomodule of $A$ which is also an ideal
of $A$. We will say that
$A$ is $H$-simple if $A$ has no nontrivial $H$-ideals.

\begin{lemma}\lelabel{2.3}
If $A$ is commutative $H$-simple, then $A^{{\rm co}H}$ is a field.
\end{lemma}

\begin{proof} Let $a$ be a nonzero element in $A^{{\rm co}H}$. Then $Aa$ is a
nonzero $H$-ideal of $A$.
But $A$ is $H$-simple and $Aa\not=0$, so $Aa=A$. Hence, we can find an
element $b$ in $A$ such that
$ba=1$.
\end{proof}

\thref{2.1} gives necessary and sufficient conditions for the projectivity of
$M\in {}_{A}\Mm^H$ as a $B$-module. These conditions might be
difficult to check, and this is why we look for sufficient conditions
that are easier.\\
$M\in {}_{A}\Mm^H$ is called coinvariantly generated if $M=AM^{{\rm co}H}$.

\begin{lemma}\lelabel{2.4}
For any $B$-module $M$, $A\otimes_{B}M$ is
coinvariantly generated. In
particular, $A$ is coinvariantly generated.
\end{lemma}

\begin{proof} 
Consider the $k$-linear map 
$$ f:\
A(A\otimes_{B}M)^{coH} \rightarrow
A\otimes_{B}M;~~f(a(u\otimes m))=au \otimes m$$ 
If $u \otimes m \in
(A\otimes_{B}M)^{coH}$, then $$u_0 \otimes m \otimes u_1=(u \otimes m)_0
\otimes (u \otimes m)_1=u
\otimes m\otimes 1$$ 
From this relation, we deduce that
\begin{eqnarray*}
&&\hspace*{-2cm}
(au \otimes m)_0 \otimes (au \otimes m)_1=a_0 u_0 \otimes m \otimes a_1
u_1\\
&=& a_0 u \otimes m
\otimes a_1=a_0(u \otimes m) \otimes a_1\\
&=&(f \otimes id_H)(a_0(u \otimes
m)_0 \otimes a_1(u
\otimes m)_1)
\end{eqnarray*}
proving that $f$ is $H$-colinear. It is clear that $f$ is $A$-linear,
and the
$k$-linear map 
$$g:\ A\otimes_{B}M \rightarrow A(A\otimes_{B}M)^{{\rm co}H}; \quad
u \otimes m \mapsto u(1
\otimes m)$$ 
is a left and right inverse of $f$.
\end{proof}

\begin{proposition}\prlabel{2.5}
Take $P\in {}_{A}\Mm^H$ and consider the following
conditions:
\begin{enumerate}
\item $A\otimes_{B}P$ is projective in ${}_{A}\Mm^H$;
\item there exists a projective coinvariantly generated object $M\in {}_{A}\Mm^H$ such
that
$M^{coH}\cong P$ as a $B$-module;
\item $P$ is projective as a $B$-module.
\end{enumerate}
Then 1) $\Rightarrow $ 2) $\Rightarrow $ 3). If the functor 
$(-)^{{\rm co}H}:\ {}_{A}\Mm^H\to{}_{B}\Mm$ is exact,
then conditions 1), 2) and 3) are equivalent.
\end{proposition}

Note that the functor $(-)^{{\rm co}H}:\ {}_{A}\Mm^H\to{}_{B}\Mm$ is exact
in the following situations:
\begin{itemize}
\item[-] $(-)^{{\rm co}H}:\ {}_{A}\Mm^H\to{}_{k}\Mm$ is exact; this is the
case if $k$ is a field and $H$
is a cosemisimple Hopf algebra \cite[Lemma 2.4.3]{Montgomery}.
\item[-]  $A$ is right $H$-coflat (see e.g. \cite[Lemma 22]{CaenepeelMZ02}); if $k$ is a
field, then this condition is equivalent
to $A$ to being injective in $\Mm^H$ (see e.g.\cite[Theorem 1]{CaenepeelMZ02}).
\end{itemize}

\begin{proof} 
1) $\Rightarrow $ 2). Let $f:\ B^{(I)}\rightarrow P$ be an
epimorphism in ${}_{B}\Mm$. Then $$1\otimes f:\
A\otimes_{B}B^{(I)}\rightarrow A\otimes_{B}P$$ is an epimorphism
in ${}_{A}\Mm^H$. Therefore $1\otimes f$ splits and we have the following
commutative diagram with exact rows:
$$\begin{diagram}
B^{(I)} & \rTo^{f}&P&\rTo&0\\
\dTo_{u_{B^{(I)}}}&&\dTo^{u_P}&&\\
(A\otimes_{B}B^{(I)})^{{\rm co}H}&\rTo^{(1\otimes f)^{{\rm co}H}}&
(A\otimes_{B}P)^{{\rm co}H}&\rTo&0
\end{diagram}$$
By \leref{1.2}, $u_{B^{(I)}}$ is an isomorphism. Therefore, $u_P$ is
surjective and, by \leref{1.3}, it is an isomorphism. Thus
$M=A\otimes_{B}P$ satisfies condition 2).\\
2) $\Rightarrow $3) Consider the map 
$$f:\ A^{(M^{{\rm co}H})}\rightarrow M;~~
f((a_m)_{m\in M^{{\rm co}H}})= \sum_{m\in M^{{\rm co}H}}a_{m}m$$
For each $m\in M^{{\rm co}H}$, we have ${\rho}(m)= m \otimes 1$, hence
$${\rho}_{M}(\sum_{m \in M^{coH}}(a_mm))=\sum_{m \in
M^{coH}}{{(a_m)}_0}m\otimes {(a_m)}_1$$
On the other hand
$$(a_m)_{m\in M^{{\rm co}H}}=\sum_{m \in M^{{\rm co}H}}( 0, 0, \cdots,
a_m, 0, \cdots, 0)$$
The linearity of ${\rho}$ implies that 
\begin{eqnarray*}
&&\hspace*{-2cm}
((a_m)_{m\in M^{{\rm co}H}})_0 \otimes (
(a_m)_{m\in
M^{coH}})_1\\
&=&\sum_{m\in M^{{\rm co}H}}(0, \cdots 0, a_m, 0, \cdots 0)_0 \otimes (0, \cdots
0, a_m, 0,  \cdots 0)_1\\
&=& \sum_{m\in M^{{\rm co}H}}(0, \cdots 0, (a_m)_0, 0, \cdots 0) \otimes (a_m)_1
\end{eqnarray*}
Thus
\begin{eqnarray*}
&&\hspace*{-2cm}
(f\otimes id_H)(((a_m)_{m\in M^{{\rm co}H}})_0 \otimes ((a_m)_{m\in
M^{{\rm co}H}})_1)\\
&=&(f\otimes
id_H)(\sum_{m\in M^{{\rm co}H}}(0, \cdots 0, (a_m)_0, 0, \cdots 0) \otimes
(a_m)_1)\\
&=&\sum{{(a_m)}_0}m\otimes
{(a_m)}_1
\end{eqnarray*}
and it follows that $f$ is $H$-colinear.  It is clear that $f$ is 
$A$-linear, so it is
a morphism in ${}_{A}\Mm^H$. $f$ is an epimorphism by the
assumption $M=AM^{{\rm co}H}$.
$f$ splits since $M$ is a projective object in $_{A}\Mm^H$,
so $M$ is a direct
summand in ${}_{A}\Mm^H$ of a direct sum of copies of $A$. By \thref{2.1},
$P$ is projective as a $B$-module.\\
3) $\Rightarrow $1)
$P$ is projective in $_{B}\Mm$, so the functor
$\Hom_B(P, -):\ {}_{B}\Mm\to {}_{k}\Mm$ is exact. On the other
hand, we have an adjoint pair of functors $(A\otimes_B (-),
(-)^{coH})$ and a natural isomorphism
$${}_A\Hom^H(A\otimes_{B}P, -) \cong \Hom_{B}(P, (-)^{{\rm co}H})$$
By assumption, the functor
$(-)^{{\rm co}H}$ is exact, and it follows that
$${}_A\Hom^H(A\otimes_{B}P, -):\ {}_A\Mm^H\to {}_{k}\Mm$$
is exact, since it is isomorphic to the composition of two exact
functors. This means that $A\otimes_{B}P$ is a projective object
in ${}_A\Mm^H$.
\end{proof}

\section {Projective and injective dimension in ${}_{A}\Mm^H$}\selabel{3}

\begin{lemma}\lelabel{3.1}
Assume that
$H$ and $A$ are commutative, and take $M, N$, $P\in {}_A\Mm^H$,
with $N$ finitely
generated projective as an $A$-module.
\begin{enumerate}
\item We have a $k$-isomorphism 
$${}_A\Hom^H( M , {}_A\Hom( N, P) )\cong
{}_A\Hom^H( M\otimes _{A}N ,P)$$
\item The functor ${}_A\Hom(N, -):\ {}_{A}\Mm^H\to {}_{A}\Mm^H$ preserves injectives.
\end{enumerate}
\end{lemma}

\begin{proof} 1) We have a natural isomorphism 
$$\phi :\ {}_A\Hom( M ,{}_A\Hom( N, P))
\rightarrow {}_A\Hom( M\otimes_{A}N , P)$$ 
$$\phi(f)(m\otimes n)=f(m)(n)$$
A standard computation shows that $f$ is $H$-colinear if and only
if $\phi(f)$ is $H$-colinear.\\
2) Let $I$ be an injective object in ${}_{A}\Mm^H$. Then the functor 
$${}_A\Hom^H( - ,I):\ {}_{A}\Mm^H\to {}_k\Mm$$
is exact. $N$ is projective as an $A$-module, so 
$(-)\otimes _{A}N$ is also exact, and it follows from 1) that
$${}_A\Hom^H( - , {}_A\Hom (N, I)):\ {}_{A}\Mm^H\to {}_{k}\Mm$$
is exact.
\end{proof}

Let $k$ be a field. Then ${}_{A}\Mm^H$ is a Grothendieck category with
enough injective
objects, and for any $M\in {}_{A}\Mm^H$ we can consider the right
derived functors ${}_A{\Ext^H}^i(M,-)$ of
${}_A{\Hom^H}(M,-):\ {}_A\Mm^H\rightarrow {}_k\Mm$.

\begin{proposition}\prlabel{3.2}
Let $k$ be a field, and assume that $H$ and $A$ are commutative. 
Take $M,N,P\in {}_{A}\Mm^H$, with $N$ finitely generated projective
as an $A$-module. Then
$${}_A{\Ext^H}^i(M,{}_A\Hom(N,P))\cong {}_A{\Ext^H}^i(M\ot_A N,P)$$
\end{proposition}

\begin{proof} By \leref{3.1}, the functors 
${}_A\Hom^H( M , {}_A\Hom( N, -) )$ and\\
${}_A\Hom^H( M\otimes _{A}N ,-)$ are isomorphic. ${}_A\Hom( N, -) $
is exact, since $N$ is projective, and it preserves injectives, by
\leref{3.1}. Therefore it preserves injective resolutions.
\end{proof}

Let ${}_A{\rm pdim}^H(-)$ and ${}_A{\rm injdim}^H(-)$ denote
respectively
the projective and injective dimension in ${}_{A}\Mm^H$.

\begin{corollary}\colabel{3.3}
With assumptions as in \prref{3.2}, we have
\begin{eqnarray}
&&{}_A{\rm pdim}^H(M\ot_A N)\leq {}_A{\rm pdim}^H(M)\\
&&{}_A{\rm injdim}^H({}_A\Hom(N,P))\leq {}_A{\rm injdim}^H(P)
\end{eqnarray}
\end{corollary}

\begin{remarks}\relabel{3.4}
1) The conclusions of \leref{3.1}, \prref{3.2} and \coref{3.3} remain
valid without the assumption that $N$ is projective, if we assume that $A$ is
semisimple.\\ 2) As a consequence of the Fundamental Theorem for Hopf modules 
(\cite[1.9.4]{Montgomery}, the results of 
\prref{3.2} and \coref{3.3} also hold without the assumption that
$N$ is projective, in the case where $A=H$.
\end{remarks}

\section {Semisimplicity of the category of relative Hopf modules}\selabel{4}
Throughout this Section, $k$ will be a field, and $A$ a $k$-algebra.
If $V$ is a finite dimensional vector space, then $A \otimes V$ is finitely generated
projective as an $A$-module. If $V$ is an $H$-comodule,
then $A\ot V\in {}_A\Mm^H$. So if $N\in {}_A\Mm^H$, then, by
\leref{1.1}, ${}_A\Hom(A\otimes V, N)\in \Mm^H$, and
${}_A\Hom(A\otimes V, N)$ and $\Hom(V, N)$ are
isomorphic as $H$-comodules.\\
Let $V$ be a finite-dimensional projective object of $\Mm^H$.
Then $A\ot V$ is a projective object of ${}_A\Mm^H$; this 
follows from the fact that
\begin{eqnarray*}
&&\hspace*{-2cm}{}_A\Hom^H(A \otimes V , N)\cong {}_A\Hom(A \otimes V , N)^{{\rm co}H}\\
&\cong & \Hom(V , N)^{{\rm co}H}\cong \Hom^H(V , N)
\end{eqnarray*}
for all $N\in {}_A\Mm^H$.

\begin{proposition}\prlabel{4.1}
$M\in {}_A\Mm^H$ is finitely
generated as an $A$-module if and only if there exist a finite-dimensional
$H$-comodule $V$ and an epimorphism of $H$-modules and $H$-comodules $\pi:\
A\otimes V\longrightarrow M$.
\end{proposition}

\begin{proof} Assume that $V$ and $\pi$ exist. Then
$A \otimes V$ is finitely
generated as $A$-module and $M$ is a quotient of $A \otimes V$ in
${}_{A}\Mm$, and therefore
finitely generated as a left $A$-module.\\
Conversely, let $M$ be finitely generated as an
$A$-module by $\{m_1,\cdots,m_n\}$. For each $i$, we can find
a finite dimensional $H$-subcomodule $W_i$ of $M$ containing $m_i$,
see \cite[5.1.1]{Montgomery}. $V=\sum W_i$
is a finite dimensional $H$-subcomodule of $M$ containing the
$m_i$ and the $k$-linear
map 
$$\pi:\ A \otimes V \rightarrow M;~~\pi(a\otimes v)=av$$
is an epimorphism of $A$-modules and $H$-comodules.
\end{proof}

Let $H^*$ be the linear dual of $H$, and let $M$ and $N$ be two right
$H$-comodules. Then $\Hom_k(M , N)$ is a left $H^*$-module, under the
action
$$(h^*f)(m)=h^*(S^{-1}(m_1)f(m_0)_1)f(m_0)$$
for all $h^*\in H^*$, $f\in \Hom_k(M , N)$ and $m\in M$.
If $M$ and $N$ are relative Hopf modules, then ${}_A\Hom(M,N)$ is
a left $H^*$-submodule of $\Hom_k(M , N)$ (see the proof of
\cite[Propositions 1.1 and 1.5]{Sweedler}.\\
A left $H^*$-module $M$ is called rational if the $H^*$-action comes
from a right $H$-coaction.

\begin{proposition}\prlabel{4.2}
Take $M,N\in {}_A\Mm^H$, and assume that $M$ is finitely generated as
an $A$-module. Then ${}_A\Hom(M, N)$ is a right $H$-comodule.
\end{proposition}

\begin{proof} Let $V$ and $\pi$ be as in \prref{4.1}. The map
$${}_A\Hom_A(\pi, N):\ {}_A\Hom(M, N)\rightarrow {}_A\Hom(A\otimes V, N)$$
is injective, and it is easy to show that it is $H^*$-linear.
${}_A\Hom(A\otimes V, N)$ is an $H$-comodule, by \leref{1.1}, so it
is a rational $H^*$-module, and - being a submodule - ${}_A\Hom(M, N)$ 
is also a rational $H^*$-module. The $H^*$-action on it is induced by
an $H$-coaction, making ${}_A\Hom(M, N)$ a right $H$-comodule.
\end{proof}

It follows from \prref{4.2} that if $A$ is commutative and if $M, N\in
{}_{A}\Mm^H$, with $M$ finitely generated as an $A$-module, then
${}_A\Hom(M , N)\in {}_{A}\Mm^H$.\\
Recall that $M\in {}_{A}\Mm^H$ is called simple if it has no proper
subobjects; a direct sum of simple objects is called semisimple.
${}_{A}\Mm^H$ is termed semisimple if every $M\in {}_{A}\Mm^H$ is
semisimple.\\
We say that ${}_{A}\Mm^H$ satisfies condition ($\dagger$) if
the functor
$${}_A\Hom(M,-):\ {}_{A}\Mm^H\to \Mm^H$$
is exact, for
every $M\in {}_{A}\Mm^H$ that is finitely generated as an $A$-module.\\
Note that ${}_{A}\Mm^H$ satisfies condition ($\dagger$) if $A$ is
semisimple; it follows from the Fundamental Theorem for Hopf modules
(\cite[1.9.4]{Montgomery}) that
${}_{H}\Mm^H$ satisfies condition ($\dagger$) if $H$
is commutative.

\begin{proposition}\prlabel{4.3}
Assume that ${}_{A}\Mm^H$ satisfies condition ($\dagger$), and that
one of the two following conditions holds:
\begin{enumerate}
\item $(-)^{{\rm co}H}: \Mm^H\to {}_k\Mm$ is exact;
\item $A$ and $H$ are commutative, and
$(-)^{{\rm co}H}: {}_A\Mm^H\to {}_B\Mm$ is exact.
\end{enumerate}
If $M\in {}_{A}\Mm^H$ is finitely
generated as an $A$-module, then it is a projective object in ${}_{A}\Mm^H$
\end{proposition}

\begin{proof}
1) The functor ${}_A\Hom^H(M,-)$ is exact since it is the composition
of the exact functors ${}_A\Hom(M,-)$ and $(-)^{{\rm co}H}$.\\
2) It follows from condition ($\dagger$) and the fact that
$H$ and $A$ are commutative that
$${}_A\Hom(M,-):\ {}_A\Mm^H\to {}_A\Mm^H$$
is exact. Now $(-)^{{\rm co}H}: {}_A\Mm^H\to {}_B\Mm$ 
and the restriction of scalars functor ${}_B\Mm\to {}_k\Mm$
are both exact, and ${}_A\Hom^H(M,-)$ is again exact, being the
composition of three exact functors.
\end{proof}

\begin{corollary}\colabel{4.4}
Let $A$ and $H$ be as in \prref{4.3}, and assume moreover that
$A$ is left noetherian. Then every $M\in {}_{A}\Mm^H$ which is
finitely generated as an $A$-module is
the direct sum of
a family of simple subobjects that are finitely generated as
$A$-modules, and consequently $M$ is a semisimple object in
${}_A\Mm^H$.
\end{corollary}

\begin{proof}
Let $N$ be a subobject of $M$ in ${}_{A}\Mm^H$. Then 
$N$ and $M/N$ are
finitely generated $A$-modules, since $A$ is left noetherian. 
It follows from \prref{4.3} that $N$ and $M/N$ are projective
objects in ${}_{A}\Mm^H$, and the exact sequence
$$0\rightarrow N\rightarrow
M \rightarrow M/N\rightarrow 0$$
splits in ${}_{A}\Mm^H$.
\end{proof}

Let $V$ be a right $H$-subcomodule of $M\in {}_{A}\Mm^H$.
It is clear that $AV$ is a subobject of $M$ in ${}_{A}\Mm^H$.

\begin{corollary}\colabel{4.5}
Let $A$ and $H$ be as in \coref{4.4}. Then ${}_{A}\Mm^H$ is
a semisimple category.
\end{corollary}

\begin{proof} 
It is well-known (see e.g. \cite[5.1.1]{Montgomery}) that every
$m\in M$ is contained in a
finite-dimensional $H$-subcomodule $V_m$ of $M$. $AV_m$
is then finitely generated as an $A$-module, and, by \coref{4.4},
the direct sum in ${}_A\Mm^H$
of a family of simple subobjects of $AV_m$ (and of $M$), all finitely generated as
$A$-modules. In particular, every $m\in M$ is contained in a simple subobject
of $M$, and this implies that $M$ is the sum of a family of simple
subobjects. Since the intersection of two simple subobjects is trivial,
it follows that this sum is a direct sum.
\end{proof}

\begin{corollary}\colabel{4.6}
Let $A$ be left noetherian semisimple and $H$
cosemisimple. Then every
object $M$ of ${}_A\Mm^H$ is a direct sum in ${}_A\Mm^H$ of a
family of simple
subobjects of $M$ finitely generated as $A$-modules. Hence $M$ is a
semisimple object in
${}_A\Mm^H$ and ${}_A\Mm^H$ is a semisimple category.
\end{corollary}

\begin{corollary}\colabel{4.7}
Let $H$ be left noetherian semisimple and cosemisimple or
commutative
noetherian, then every object $M$ of ${}_H\Mm^H$ is a direct sum in
${}_H\Mm^H$ of a
family of simple subobjects of $M$ finitely generated as $H$-modules. Hence
$M$ is a semisimple
object in ${}_H\Mm^H$ and ${}_H\Mm^H$ is a semisimple category.
\end{corollary}

\end{document}